\documentclass{gtart}

\def\ifplaintex{\expandafter\ifx\csname documentclass\endcsname\relax}

\expandafter\ifx\csname epsfbox\endcsname\relax\input epsf\fi

\ifplaintex 
\else
\fi

\def\gt{{\mathsurround=0pt\it $\cal G\mskip-2mu$eometry \&\ 
$\cal T\!\!$opology}}        

\def\gtp{{\mathsurround=0pt\it $\cal G\mskip-2mu$eometry \&\ 
$\cal T\!\!$opology $\cal P\!$ublications}}  

\def\lognumber#1{\def\thelognumber{#1}}
\def\volumenumber#1{\def\thevolumenumber{#1}}
\def\papernumber#1{\def\thepapernumber{#1}}
\def\volumeyear#1{\def\thevolumeyear{#1}}

\def\pagenumbers#1#2{\def\startpage{#1}\def\finishpage{#2}}
\def\published#1{\def\publishdate{#1}}
\def\proposed#1{\def\theproposer{#1}}
\def\seconded#1{\def\theseconders{#1}}
\def\received#1{\def\receiveddate{#1}}

\def\accepted#1{\def\accepteddate{#1}}
\def\asciititle#1{\def\theasciititle{#1}}

\long\def\asciiabstract#1{\long\def\theasciiabstract{#1}}
\def\asciikeywords#1{\def\theasciikeywords{#1}}

\let\\\par\let\thevolumenumber\relax\let\thepapernumber\relax
\let\thevolumeyear\relax\let\thesamplenumber\relax\let\startpage\relax
\let\finishpage\relax\let\publishdate\relax\let\receiveddate\relax
\let\reviseddate\relax\let\accepteddate\relax\let\theasciititle\relax
\let\theasciiauthors\relax
\let\theasciiabstract\relax\let\theasciikeywords\relax
\let\theasciiemail\relax\let\theshortauthors\relax\let\theshorttitle\relax

\long\def\maketitlep{   

\count0=\startpage

\gt\hfill      
\hbox to 77pt{\vbox to 0pt{\vglue -15pt\epsfbox{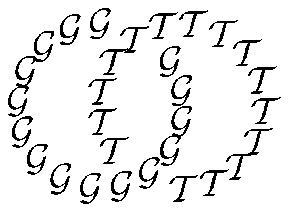}\vss}\hss}
\break
{\small\ifx\thesamplenumber\relax 
Volume \else Sample
\fi\thevolumenumber\ (\thevolumeyear)
\startpage--\finishpage\nl
Published: \publishdate}
\vglue 0.5truein plus 0.4fil minus 0.1truein

{\parskip=0pt\leftskip 0pt plus 1fil\def\\{\par\smallskip}{\ifplaintex\large
\else\Large\fi\bf\thetitle}\par\medskip}   

\vglue 0pt plus 0.1fil 

{\parskip=0pt\leftskip 0pt plus 1fil\def\\{\par}{\sc\theauthors}
\par\medskip}

\vglue 0pt plus 0.1fil 

{\small\parskip=0pt\let\newline\\
{\leftskip 0pt plus 1fil\def\\{\par}{\sl\theaddress}\par}
\expandafter\ifx\theemail\relax    
\relax\else\vglue 5pt plus 0.02fil minus 2pt\def\\{\stdspace{\rm 
and}\stdspace} 
\cl{Email:\stdspace\tt\theemail}\fi
\ifx\theurl\relax                  
\relax\else\vglue 5pt plus 0.02fil minus 2pt\def\\{\stdspace{\rm 
and}\stdspace}
\cl{URL:\stdspace\tt\theurl}\fi\par}

\vglue 7pt plus 0.3fil minus 3pt

{\bf Abstract}
\vglue 5pt plus 0.1fil minus 2pt

\theabstract

\vglue 7pt plus 0.3fil minus 3pt

{\bf AMS Classification numbers}\quad Primary:\quad \theprimaryclass

Secondary:\quad \thesecondaryclass

\vglue 5pt plus 0.3fil minus 2pt

{\bf Keywords:}\quad \thekeywords

\vglue 10pt plus 0.5fil minus 5pt

{\small  Proposed: \theproposer\hfill Received: \receiveddate\nl
Seconded: \theseconders\hfill 
\ifx\reviseddate\relax                         
Accepted: \accepteddate                        
\else
Revised: \reviseddate                          
\fi}
\eject
}       

\let\maketitlepage\maketitlep
\let\maketitle\maketitlepage

\font\phead=cmsl9 scaled 950
\font=cmsl9 scaled 1050
\font\pnum=cmbx10 scaled 913
\font\lnum=cmbx10 
\font\pfoot=cmsl9 scaled 950
\font=cmsl9 scaled 1050
\ifplaintex
\headline{\vbox to 0pt{\vskip -4.5mm\line{\small\phead\ifnum
\count0=\startpage ISSN 1364-0380 (on line)
1465-3060 (printed) \hfill {\pnum\folio}\else\ifodd\count0\def\\{ }%
\ifx\theshorttitle\relax\thetitle\else\theshorttitle\fi\hfill{\pnum\folio}
\else\def\\{ and }{\pnum\folio}\hfill\ifx\theshortauthors\relax\theauthors
\else\theshortauthors\fi\fi\fi}\vss}}
\footline{\vbox to 0pt{\vglue 0mm\line{\small\pfoot\ifnum\count0=\startpage
\copyright\ \gtp\hfill\else
\gt, Volume \thevolumenumber\ (\thevolumeyear)\hfill\fi}\vss
}}
\else
\makeatletter
\def\@oddhead{{\small\lhead\ifnum\count0=\startpage ISSN 1364-0380 (on line)
1465-3060 (printed) \hfill {\lnum\number\count0}\else\ifodd\count0
\def\\{ }\ifx\theshorttitle\relax \thetitle \else\theshorttitle\fi\hfill
{\lnum\number\count0}\else\def\\{ and }{\lnum\number\count0}
\hfill\ifx\theshortauthors\relax 
\theauthors\else\theshortauthors\fi\fi\fi}}\def\@evenhead{@oddhead}
\def\@oddfoot{\small\lfoot\ifnum\count0=\startpage\copyright\ \gtp\hfill\else
\gt, Volume \thevolumenumber\ (\thevolumeyear)\hfill\fi}
\def\@evenfoot{@oddfoot}
\makeatother
\fi

\newwrite\gtoutfile
\long\gdef\makeheadfile{  
{\def\\{, }\def\s{ }
\immediate\openout\gtoutfile head.xxx
\immediate\write\gtoutfile{To: math@arxiv.org}
\immediate\write\gtoutfile{Subject: put OR rep NNNNN:pppp}
\immediate\write\gtoutfile{--text follows this line--}
\immediate\write\gtoutfile{Proxy-for: \ifx\theasciiauthors\relax
\theauthors\else\theasciiauthors\fi\s<\ifx\theasciiemail\relax\theemail\else\theasciiemail\fi>}
\immediate\write\gtoutfile{\noexpand\\}
\immediate\write\gtoutfile{Authors: \ifx\theasciiauthors\relax
\theauthors\else\theasciiauthors\fi}
{\def\\{ }\immediate\write\gtoutfile{Title: \ifx\theasciititle\relax
\thetitle\else\theasciititle\fi}}
\immediate\write\gtoutfile{Subj-class: GT or GR or SG or ...}
\immediate\write\gtoutfile{MSC-class: \theprimaryclass\ifx\thesecondaryclass\relax\else, \thesecondaryclass\fi}
\immediate\write\gtoutfile{Journal-ref: Geom. Topol. \thevolumenumber\s
(\thevolumeyear) \startpage-\finishpage}
\immediate\write\gtoutfile{Comments: Published in Geometry and Topology at}
\immediate\write\gtoutfile{    http://www.maths.warwick.ac.uk/gt/GTVol\thevolumenumber/paper\thepapernumber.abs.html}
\immediate\write\gtoutfile{\noexpand\\}
\immediate\write\gtoutfile{}
\ifx\theasciiabstract\relax
\immediate\write\gtoutfile{\theabstract}\else
\immediate\write\gtoutfile{\theasciiabstract}\fi
\immediate\write\gtoutfile{}
\immediate\write\gtoutfile{\noexpand\\}
\immediate\write\gtoutfile{}
\immediate\closeout\gtoutfile}}  

\def\maketitlepage{\maketitlep\makeheadfile}
\let\maketitle\maketitlepage

\lognumber{240}

\volumenumber{6}
\papernumber{15}
\volumeyear{2002}
\pagenumbers{409}{424}
\received{27 February 2002}
\accepted{22 August 2002}
\published{15 September 2002}
\proposed{Cameron Gordon}

\seconded{Joan Birman, Walter Neumann}

\usepackage{amssymb,graphicx,amsmath,diagrams}

\begingroup\makeatletter\ifx\SetFigFont\undefined%
\gdef\SetFigFont#1#2#3#4#5{%
  \reset@font\fontsize{#1}{#2pt}%
  \fontfamily{#3}\fontseries{#4}\fontshape{#5}%
  \selectfont}%
\fi\endgroup%

\newtheorem{theorem}{Theorem}[section]
\newtheorem*{maintheorem}{Theorem \ref{mainthm}}
\newtheorem*{secondprop}{Proposition \ref{lcs}}
\newtheorem{corollary}[theorem]{Corollary}
\newtheorem{proposition}[theorem]{Proposition}
\theoremstyle{remark}
\newtheorem*{acknowledgement}{Acknowledgements}
\newtheorem{example}[theorem]{Example}
\newtheorem{remark}[theorem]{Remark}
\newtheorem{question}[theorem]{Question}

\newcommand{\rank}{{\rm rank}}

\begin{document}

\title{On the Cut Number of a 3--manifold}
\asciititle{On the Cut Number of a 3-manifold}
\author{Shelly L Harvey}
\address{Department of Mathematics\\University of
California at San Diego\\La Jolla, CA 92093-0112, USA}
\email{sharvey@math.ucsd.edu}
\url{http://math.ucsd.edu/\char'176sharvey}
\asciikeywords{3-manifold, fundamental group, corank, Alexander
module, virtual betti number, free group}

\begin{abstract}
The question was raised as to whether the cut number of a 3--manifold $X$
is bounded from below by $\frac{1}{3}\beta _{1}\left( X\right) $.
We show that the answer to this question is ``no.''  For each $m\geq 1$, we
construct explicit examples of closed 3--manifolds $X$ with $\beta
_{1}\left( X\right) =m$ and cut number 1.  That is, $\pi_1\left(X\right)$
cannot map onto any non-abelian free group.  Moreover, we show that these
examples can be assumed to be hyperbolic.
\end{abstract}

\asciiabstract{The question was raised as to whether the cut number
of a 3-manifold X is bounded from below by 1/3 beta_1(X).  We show
that the answer to this question is `no.'  For each m>0, we construct
explicit examples of closed 3-manifolds X with beta_1(X)=m and cut
number 1.  That is, pi_1(X) cannot map onto any non-abelian free
group.  Moreover, we show that these examples can be assumed to be
hyperbolic.}

\maketitlepage

\section{\label{intro}Introduction}

Let $X$ be a closed, orientable $n$--manifold.  The \emph{cut number of }$X$, $
c\left( X\right) $, is defined to be the maximal number of components of a
closed, 2--sided, orientable hypersurface $F\subset X$ such that $X-F$ is
connected.  Hence, for any $n\leq c\left( X\right) $, we can construct a
map $f\co X\rightarrow \bigvee_{i=1}^{n}S^{1}$ such that the induced map on $
\pi _{1}$ is surjective.  That is, there exists a surjective map $f_{\ast
}\co \pi _{1}\left( X\right) \twoheadrightarrow F\left( c\right) $, where $
F\left( c\right) $ is the free group with $c=c\left( X\right) $ generators.
Conversely, if we have any epimorphism $\phi \co \pi _{1}\left( X\right)
\twoheadrightarrow F\left( n\right) $, then we can find a map $
f\co X\rightarrow \bigvee_{i=1}^{n}S^{1}$ such that $f_{\ast }=\phi $. After
making the $f$ transverse to a non-wedge point $x_{i}$ on each $S^{1}$, $
f^{-1}\left( X\right) $ will give $n$ disjoint surfaces $F=\cup F_{i}$ with $
X-F$ connected. Hence one has the following elementary group-theoretic
characterization of $c\left( X\right)$.

\begin{proposition}
\label{group}$c\left( X\right) $ is the maximal $n$ such that there is an
epimorphism $\phi \co \pi _{1}\left( X\right) \twoheadrightarrow F\left(
n\right) $ onto the free group with $n$ generators.
\end{proposition}

\begin{example}
Let $X=S^{1}\times S^{1}\times S^{1}$ be the 3--torus. Since $\pi
_{1}\left( X\right) =\mathbb{Z}^{3}$ is abelian, $c\left( X\right) =1$.
\end{example}

Using Proposition \ref{group}, we show that the cut number is additive under
connected sum.

\begin{proposition}\label{connectsum}
If $X=X_{1}\#X_{2}$ is the connected sum of $X_{1}$ and $X_{2}$ then
\begin{equation*}
c\left( X\right) =c\left( X_{1}\right) +c\left( X_{2}\right) \text{.}
\end{equation*}
\end{proposition}

\begin{proof}
Let $G_{i}=\pi _{1}\left( X_{i}\right) $ for $i=1,2$ and $G=\pi _{1}\left(
X\right) \cong G_{1}\ast G_{2}$. It is clear that $G$ maps surjectively onto $F\left( c\left( X_{1}\right) \right) \ast F\left( c\left( X_{2}\right) \right) \cong F\left( c\left( X_{1}+X_{2}\right) \right) $. 
Therefore
$c\left(  X \right) \geq c\left( X_{1}\right) +c\left( X_{2}\right)$.

Now suppose that there exists a map $\phi \co G\twoheadrightarrow F\left( n\right) $.
Let $\phi _{i}\co G_{i}\rightarrow F\left( n\right) $ be the composition $G_{i}\rightarrow G_{1}\ast G_{2}\overset{\cong }{\rightarrow }G\overset{\phi }{\twoheadrightarrow }F\left( n\right) $.
 Since $\phi $ is surjective and $G\cong G_{1}\ast G_{2}$,
$\text{Im}\left( \phi
_{1}\right) $ and $\text{Im}\left( \phi _{2}\right) $ generate $F\left(
n\right) $.
 Morever, $\text{Im}\left( \phi _{i}\right) $ is a subgroup of a free
group, hence is free of rank less than or
equal to $c\left( X_{i}\right) $. It follows that $n\leq c\left( X_{1}\right) +c\left( X_{2}\right) $.  In particular, when $n$ is maximal we have
$c\left( X\right) =n\leq c\left( X_{1}\right) +c\left( X_{2}\right)$.
\end{proof}

In this paper, we will only consider 3--manifolds with $\beta _{1}\left(
X\right) \geq 1$. Consider the surjective map $\pi _{1}\left( X\right)
\twoheadrightarrow H_{1}\left( X\right) /\left\{ \mathbb{Z}\text{--torsion}
\right\} \cong \mathbb{Z}^{\beta _{1}\left( X\right) }$ . Since $\beta
_{1}\left( X\right) \geq 1$, we can find a surjective map from $\mathbb{Z}
^{\beta _{1}\left( X\right) }$ onto $\mathbb{Z}$. It follows from
Proposition \ref{group} that $c\left( X\right) \geq 1$. Moreover, every
map $\phi \co \pi _{1}\left( X\right) \twoheadrightarrow F\left( n\right) $
gives rise to an epimorphism $\overline{\phi}\co H_{1}\left( X\right)
\twoheadrightarrow H_{1}\left(\bigvee_{i=1}^{n}S^{1}\right) \cong
\mathbb{Z}^{n}$ It follows that $\beta_{1}\left( X\right) \geq n$ which
gives us the well known result:
\begin{equation}
1\leq c\left( X\right) \leq \beta _{1}\left( X\right) \text{.}
\end{equation}

It has recently been asked whether a (non-trivial) lower bound exists for the
cut number.  We make the following observations.

\begin{remark}
If $S$ is a closed, orientable surface then $c\left( S\right)=\frac{1}{2}\beta
_{1}\left( S\right) $.
\end{remark}

\begin{remark}
If $X$ has solvable fundamental group then $c\left( X\right) =1$
and $\beta _{1}(X)$ $\leq 3$.
\end{remark}

\begin{remark}
Both $c$ and $\beta _{1}$ are additive under connected sum
(Proposition \ref
{connectsum}).
\end{remark}

Therefore it is natural to ask the following question first asked
by A Sikora and T Kerler.  This question was motivated by
certain results and conjectures on the divisibility of quantum
3--manifold invariants by P Gilmer--T Kerler
\cite{gilmer_kerler:cutnumber} and T Cochran--P Melvin
\cite{cochran_melvin:quantum}.
\begin{question}
\label{cutquest}Is $c\left( X\right) \geq \frac{1}{3}\beta _{1}\left(
X\right) $ for all closed, orientable 3--manifolds $X$?
\end{question}
We show that the answer to this question is ``as far from yes as possible.''
In fact, we show that for each $m\geq 1$ there exists a closed,
\emph{hyperbolic} 3--manifold with $\beta _{1}\left(
X\right) =m$ and $c\left( X\right) =1$. We actually prove a stronger
statement.

\begin{maintheorem}
For each $m\geq 1$ there exist closed 3--manifolds $X$ with $
\beta _{1}\left( X\right)$ $=m$ such that for any infinite cyclic cover $
X_{\phi }\rightarrow X$, $\rank_{\mathbb{Z}\left[ t^{\pm 1}\right]
}H_{1}\left( X_{\phi }\right) =0$.
\end{maintheorem}

We note the condition stated in the Theorem \ref{mainthm} is
especially interesting because of the following theorem of J
Howie \cite{howie:free_subgroups}. Recall that a group $G$ is
\emph{large} if some subgroup of finite index has a non-abelian
free homomorphic image.  Howie shows that if $G$ has an infinite
cyclic cover whose rank is at least $1$ then $G$ is large.
\begin{theorem}[Howie \protect\cite{howie:free_subgroups}]
Suppose that $\widetilde{K}$ is a connected regular covering complex of a
finite 2--complex $K$, with nontrivial free abelian covering transformation
group $A$.  Suppose also that $H_{1}\left( \widetilde{K};\mathbb{Q}\right)
$ has a free $\mathbb{Q}\left[ A\right] $--submodule of rank at least $1$.
Then $G=\pi _{1}\left( K\right) $ is large.
\end{theorem}

Using the proof of Theorem \ref{mainthm} we show that
the fundamental group of the aforementioned 3--manifolds cannot map onto
$F/F_{4}$ where $F$ is the free group with $2$ generators and $F_{4}$ is
the $4^{th}$ term of the lower central series of $F$.

\begin{secondprop}Let $X$ be as in Theorem
\ref{mainthm}, $G=\pi_1\left(X\right)$ and $F$ be the
free group on $2$ generators.  There is no epimorphism from $G$ onto
$F/F_4$. \end{secondprop}

Independently, A Sikora has recently shown that the cut number of
a ``generic'' 3--manifold is at most 2 \cite{sikora:cutnumber}.
Also, C Leininger and A Reid have constructed specific examples
of genus 2 surface bundles $X$ satisfying (i)
$\beta _{1}\left( X\right) =5$ and $c\left( X\right) =1$ and
(ii) $\beta _{1}\left( X\right) =7$ and $ c\left(
X\right) =2$ \cite{leininger_reid:corank}.

\begin{acknowledgement}
I became interested in the question as to whether the cut number of a
3--manifold was bounded below by one-third the first betti number after
hearing it asked by A Sikora at a problem session of the 2001 Georgia
Topology Conference.  The
question was also posed in a talk by T Kerler at the 2001
Lehigh Geometry and Topology Conference.  The author was supported by NSF
DMS-0104275 as well as by the Bob E and Lore Merten Watt Fellowship.
\end{acknowledgement}

\section{Relative Cut Number}

Let $\phi $ be a primitive class in $H^{1}\left(
X;\mathbb{Z}\right) $. Since $H^{1}\left( X;\mathbb{Z}\right)
\cong \operatorname{Hom}\left( \pi _{1}\left( X\right) ,\mathbb{Z}\right)
$, we
can assume $\phi $ is a surjective homomorphism, $\phi \co \pi
_{1}\left( X\right) \twoheadrightarrow \mathbb{Z}$. Since $X$ is
an orientable 3--manifold, every element in $H_{2}\left( X;
\mathbb{Z}\right) $ can be represented by an embedded, oriented,
2--sided surface \cite[Lemma 1]{thurston:norm}. Therefore, if $\phi
\in H^{1}\left( X;\mathbb{Z} \right) \cong H_{2}\left(
X;\mathbb{Z}\right) $ there exists a surface (not unique) dual to
$\phi $.  The \emph{cut number of }$X$\emph{\ relative to }$ \phi
$, $c\left( X,\phi \right) $, is defined as the maximal number of
components of a closed, 2--sided, oriented surface $F\subset X$
such that $ X-F $ is connected and one of the components of $F$ is
dual to $\phi $.  In the above definition, we could have required
that ``any number'' of components of $F$ be dual to $\phi $ as
opposed to just ``one.''  We remark that since $X-F$ is connected,
these two conditions are equivalent. Similar to $c\left(
X\right) $, we can describe $c\left( X,\phi \right) $ group
theoretically.

\begin{proposition}
$c\left( X,\phi \right) $ is the maximal $n$ such that there is an
epimorphism $\psi \co \pi _{1}\left( X\right) \twoheadrightarrow F\left(
n\right) $ onto the free group with $n$ generators that factors through $
\phi $ (see diagram on next page).

\cl{$\begin{diagram} \pi_1\left(X\right) & \rTo^{\phi} & \mathbb{Z} \\ \dTo^\psi
& \ruTo \\ F \left( n \right) \end{diagram}$}
\end{proposition}

It follows immediately from the definitions that $c\left( X,\phi \right)
\leq c\left( X\right) $ for all primitive $\phi $.  Now let $F$ be any
surface with $c\left( X\right) $ components and let $\phi $ be dual to one
of the components, then $c\left( X,\phi \right) =c\left( X\right) $.  Hence
\begin{equation}
c\left( X\right) =\max \left\{ c\left( X,\phi \right) \mid\phi \text{ is a
primitive element of }H^{1}\left( X;\mathbb{Z}\right) \right\} \text{.}
\label{cutviachi}
\end{equation}
In particular, if $c\left( X,\phi \right) =1$ for all $\phi $ then $c\left(
X\right) =1$.

We wish to find sufficient conditions for $c\left( X,\phi \right)
=1$.  In \cite[page 44]{kerler:hom_tqft}, T Kerler develops a
skein theoretic algorithm to compute the one-variable Alexander
polynomial $\Delta _{X,\phi }$ from a surgery presentation of $X$.
As a result, he shows that if $c\left( X,\phi \right) \geq 2$ then
the Frohman--Nicas TQFT evaluated on the cut cobordism is zero,
implying that $\Delta _{X,\phi }=0$.  Using the fact that
$\mathbb{Q}\left[ t^{\pm 1}\right] $ is a principal ideal domain
one can prove that $\Delta _{X,\phi }=0$ is equivalent to
$\rank_{\mathbb{Z}\left[ t^{\pm 1}\right] }H_{1}\left( X_{\phi
}\right) \geq 1$.  We give an elementary proof of the equivalent
statement of Kerler's.

\begin{proposition}\label{cutandinf}
If $c\left( X,\phi \right) \geq 2$ then $\rank_{\mathbb{Z}\left[ t^{\pm 1}
\right] }H_{1}\left( X_{\phi }\right) \geq 1$.
\end{proposition}

\begin{proof}
Suppose $c\left( X,\phi \right)\geq 2$ then there is a surjective map
$\psi \co \pi _{1}\left( X\right) \twoheadrightarrow F\left(
n\right) $ that factors through $\phi $ with $n\geq 2$. 
Let $\overline{\phi }\co F\left( n\right) \twoheadrightarrow \mathbb{Z}$ be the
homomorphism such that $\phi =\overline{\phi }\circ \psi $.  $\phi $
surjective implies that $\psi_{|\ker \phi }\co \ker
\phi \twoheadrightarrow \ker \overline{\phi } $ is
surjective.  Writing $\mathbb{Z}$ as the multiplicative group generated by $
t$, we can consider $\frac{\ker \phi }{\left[ \ker
\phi ,\ker  \phi  \right] }$ and $\frac{\ker
\overline{\phi } }{\left[ \ker \overline{\phi } ,\ker
 \overline{\phi } \right] }$ as modules over
$\mathbb{Z}\left[t^{\pm 1}\right] $. Here, the t acts by conjugating by an element that maps to t by
$\phi$ or $\overline{\phi }$.   Moreover, $\psi_{
|\ker \phi }\co \frac{
\ker \phi }{\left[ \ker  \phi ,\ker \phi
 \right] }\twoheadrightarrow \frac{\ker \overline{\phi }
 }{\left[ \ker \overline{\phi } ,\ker \overline{
\phi } \right] }$ is surjective hence
\begin{equation*}
\rank_{\mathbb{Z}\left[ t^{\pm 1}\right] }\left( \frac{\ker \phi
 }{\left[ \ker  \phi  ,\ker  \phi \right] }
\right) \geq \rank_{\mathbb{Z}\left[ t^{\pm 1}\right] }\left( \frac{\ker
 \overline{\phi } }{\left[ \ker  \overline{\phi }
,\ker \overline{\phi } \right] }\right)=n-1 \text{.}
\end{equation*}
Since $n\geq 2$, $\rank_{\mathbb{Z}\left[ t^{\pm 1}\right] }
H_1\left(X_\phi\right)=\rank_{\mathbb{Z}\left[
t^{\pm 1}\right] }\left( \frac{\ker
 \phi  }{\left[ \ker  \phi
,\ker  \phi  \right] }\right) \geq 1$.
\end{proof}

\begin{corollary}\label{ontoFmodF2}
If $\pi_1\left(X\right) \twoheadrightarrow F/F''$ where $F$ is a free group of
rank 2 then there exists a $\phi \co  \pi_1\left(X\right) \twoheadrightarrow
\mathbb{Z}$ such that
$\rank_{\mathbb{Z}\left[ t^{\pm
1}\right] }H_{1}\left( X_{\phi }\right) \geq 1$.
\end{corollary}\eject
\begin{proof}This follows immediately from the proof of Proposition
\ref{cutandinf} after noticing that $F'' \subset \left[ \ker
\left( \overline{\phi} \right),\ker \left( \overline{\phi} \right)
\right]$ and $\operatorname{Hom}\left(F/F'',\mathbb{Z}\right)\cong
\operatorname{Hom}\left(F,\mathbb{Z}\right)$.
\end{proof}

\section{The Examples}

We construct closed 3--manifolds all of whose infinite cyclic covers have
first homology that is $\mathbb{Z}\left[ t^{\pm 1}\right] $--torsion.  The
3--manifolds we
consider are 0--surgery on an $m$--component link that is obtained from the
trivial link by tying a Whitehead link interaction between each two
components.

\begin{theorem}
\label{mainthm}For each $m\geq 1$ there exist closed
3--manifolds $X$ with $\beta _{1}\left( X\right)$ $=m$ such that for any
infinite cyclic cover $X_{\phi }\rightarrow X$, $\rank_{\mathbb{Z}\left[
t^{\pm 1}\right] }H_{1}\left( X_{\phi }\right) =0$.
\end{theorem}

It follows from Propostion \ref{cutandinf} that the cut number of the
manifolds in Theorem \ref{mainthm} is 1.  In fact,
Corollary \ref{ontoFmodF2} implies that
$\pi_1\left(X\right)$ does not map onto $F/F''$ where $F$ is a free group of
rank 2.  Moreover, the proof of this theorem shows that
$\pi_1\left(X\right)$ does not even map onto $F/F_4$ where $F_n$ is the
$n^{th}$ term of the lower central series of $F$ (see Proposition
\ref{lcs}).

By a theorem of Ruberman \cite{ruberman:seifert}, we can assume
that the manifolds with cut number $1$ are hyperbolic.

\begin{corollary}
For each $m\geq 1$ there exist closed, orientable, hyperbolic 3--manifolds $Y
$ with $\beta _{1}\left( Y\right) =m$ such that for any infinite cyclic
cover $Y_{\phi }\rightarrow Y$, $\rank_{\mathbb{Z}\left[ t^{\pm 1}\right]
}H_{1}\left( Y_{\phi }\right) =0$.
\end{corollary}

\begin{proof}
Let $X$ be one of the 3--manifolds in Theorem \ref{mainthm}.  By
\cite[Theorem 2.6]{ruberman:seifert}, there exists a degree one
map $f\co Y\rightarrow X$ where $Y$ is hyperbolic and $f_{\ast }$ is
an isomorphism on $H_{\ast }$.  Denote by $ G=\pi _{1}\left(
X\right) $ and $P=\pi _{1}\left( Y\right) $.  It is then
well-known that $f$ is surjective on $\pi _{1}$.  It follows from
Stalling's theorem \cite[page 170]{stallings:centralseries} that the
kernel of $f_{\ast }$ is $P_{\omega }\equiv \cap P_{n}$.  Now,
suppose $\phi \co P\overset{f_{\ast }}{ \twoheadrightarrow
}G\overset{\overline{\phi }}{\twoheadrightarrow }\mathbb{Z }$
defines an infinite cyclic cover of $Y$.  Then $H_{1}\left(
Y_{\phi }\right) \twoheadrightarrow H_{1}\left( X_{\overline{\phi
}}\right) $ has kernel $P_{\omega }/\left[ \ker \phi ,\ker \phi
\right] $.  To show that $ \rank_{\mathbb{Z}\left[ t^{\pm 1}\right]
}H_{1}\left( Y_{\phi }\right) =0$ it suffices to show that
$P_{\omega }$ vanishes under the map $H_{1}\left( Y_{\phi }\right)
\rightarrow H_{1}\left( Y_{\phi }\right) \otimes _{\mathbb{Z
}\left[ t^{\pm 1}\right] }\mathbb{Q}\left[ t^{\pm 1}\right]
\rightarrow H_{1}\left( Y_{\phi }\right) \otimes
_{\mathbb{Z}\left[ t^{\pm 1}\right] } \mathbb{Q}\left( t\right) $
since then $\rank_{\mathbb{Z}\left[ t^{\pm 1} \right] }H_{1}\left(
Y_{\phi }\right) =\rank_{\mathbb{Z}\left[ t^{\pm 1} \right]
}H_{1}\left( X_{\overline{\phi }}\right) =0$.

Note that $H_{1}\left( Y_{\phi }\right) \otimes _{\mathbb{Z}\left[ t^{\pm 1}
\right] }\mathbb{Q}\left[ t^{\pm 1}\right] \cong \bigoplus_{i=1}^{n}\mathbb{Q
}\left[ t^{\pm 1}\right] \oplus T$ where $T$ is a $\mathbb{Q}\left[ t^{\pm 1}
\right] $ torsion module. Moreover, $P_{n}$ is generated by elements of
the form $\gamma =\left[ p_{1}\left[ p_{2}\left[ p_{3},\ldots \left[
p_{n-2},\alpha \right] \right] \right] \right] $ where $\alpha \in
P_{2}\subseteq \ker \phi $.  Therefore
\begin{equation*}\left[ \gamma \right] =\left( \phi
\left( p_{i}\right) -1\right) \cdots \left( \phi \left( p_{n-2}\right)
-1\right) \left[ \alpha \right]
\end{equation*}
 in $H_{1}\left( Y_{\phi }\right) $ which implies that $P_{n}\subseteq
J^{n-2}\left( H_{1}\left( Y_{\phi }\right) \right)
$ for $n\geq 2$
where $J$ is the augmentation ideal of
$\mathbb{Z}\left[
t^{\pm 1}\right]$.  It follows that any element of $P_{\omega }$
considered
as an element of $H_{1}\left( Y_{\phi }\right) \otimes _{\mathbb{Z}\left[
t^{\pm 1}\right] }\mathbb{Q}\left[ t^{\pm 1}\right] $ is infinitely
divisible by $t-1$ and hence is torsion.
\end{proof}

\begin{figure}[ht!]
\begin{center}
\begin{picture}(0,0)%
\includegraphics{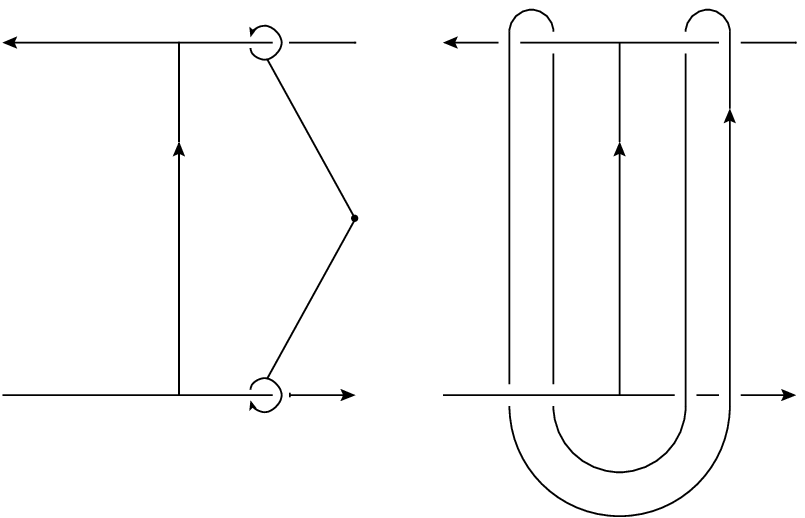}%
\end{picture}%
\setlength{\unitlength}{3947sp}%
\begin{picture}(3832,2449)(328,-1702)
\small
\put(948,-108){$\alpha_{ij}$}%

\put(400,632){$L_i$}%

\put(1874,-1300){$L_j$}%

\put(1715,632){$x_i$}%

\put(1504,-1300){$x_j$}%

\put(2500,632){$L_i$}%

\put(3994,-1300){$L_j$}%

\put(3050,-108){$\alpha_{ij}$}%

\put(3884,209){$\gamma_{ij}$}%

\end{picture}
\nocolon\caption{}
\label{figure1}
\end{center}
\end{figure}

\begin{figure}[ht!]
\vspace{-0.8in}
\begin{center}
\setlength{\unitlength}{0.01in}
\begin{picture}(500,500)(0,0)
\put(75,0){
\includegraphics{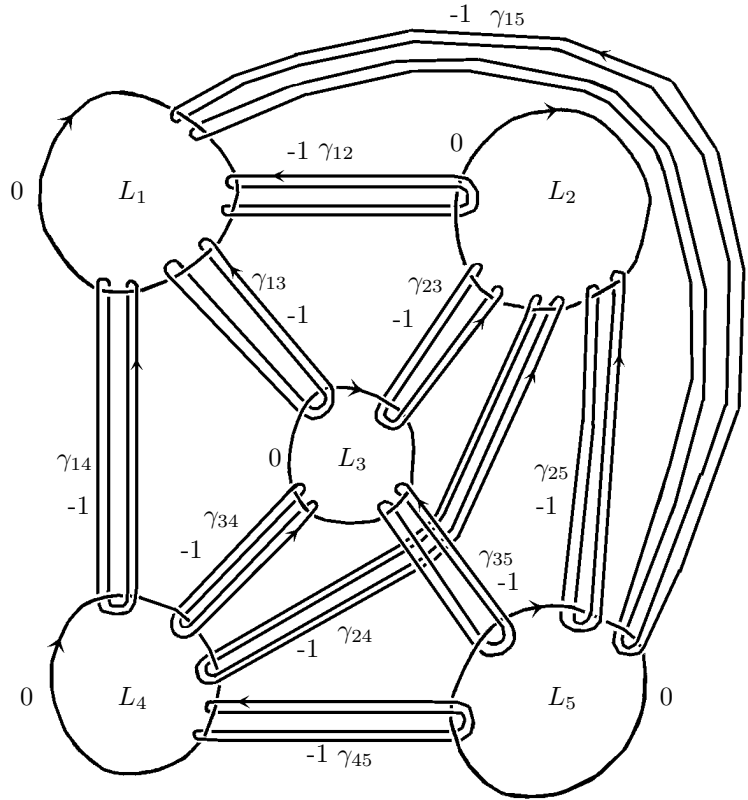}
}

\small
\put(240,175){$L_3$}
\put(205,175){0}

\put(350,50){$L_5$}%
\put(410,50){0}

\put(125,50){$L_4$}%
\put(75,50){0}

\put(125,315){$L_1$}%
\put(70,315){0}

\put(350,315){$L_2$}%
\put(300,340){0}

\put(343,170){$\gamma_{25}$}
\put(345,150){-1}

\put(320,408){$\gamma_{15}$}
\put(300,408){-1}

\put(195,270){$\gamma_{13}$}
\put(215,250){-1}

\put(277,268){$\gamma_{23}$}
\put(270,248){-1}

\put(314,125){$\gamma_{35}$}
\put(325,110){-1}

\put(240,85){$\gamma_{24}$}
\put(220,75){-1}

\put(170,145){$\gamma_{34}$}
\put(159,127){-1}

\put(230,338){$\gamma_{12}$}
\put(215,336){-1}

\put(93,175){$\gamma_{14}$}
\put(100,150){-1}

\put(240,20){$\gamma_{45}$}
\put(225,20){-1}

\end{picture}
\caption{The surgered manifold $X$ when $m=5$}
\label{figure2}
\end{center}
\end{figure}

\proof[Proof of Theorem \ref{mainthm}]
Let $L=\sqcup L_{i}$ be the oriented trivial link with $m$
components in $ S^{3}$ and $\sqcup D_{i}$ be oriented disjoint
disks with $\partial D_{i}=L_{i}$.  The fundamental group of
$S^{3}-L$ is freely generated by $ \left\{ x_{i}\right\} $ where
$x_{i}$ is a meridian curve of $L_{i}$ which intersects $D_{i}$
exactly once and $D_{i}\cdot x_{i}=1$.  For all $i,j$ with $1\leq
i<j\leq m$ let $\alpha _{ij}\co I\rightarrow S^{3}$ be oriented
disjointly embedded arcs such that $\alpha _{ij}\left( 0\right)
\in L_{i}$ and $\alpha _{ij}\left( 1\right) \in L_{j}$ and $\alpha
_{ij}\left( I\right) $ does not intersect $\sqcup D_{i}$.  For
each arc $\alpha _{ij}$, let $ \gamma _{ij}$ be the curve embedded
in a small neighborhood of $\alpha _{ij}$ representing the class
$\left[ x_{i},x_{j}\right] $ as in Figure \ref {figure1}.
Let $X$ be the 3--manifold obtained performing 0--framed Dehn
surgery on $L$ and $-1$--framed Dehn surgery on each $\gamma =\sqcup
\gamma _{ij}$.  See Figure \ref{figure2}
for an example of $X$ when $m=5$.

Denote by $X_{0}$, the manifold obtained by performing 0--framed Dehn surgery
on $L$.  Let $W$ be the 4--manifold obtained by adding a 2--handle to $
X_{0}\times I$ along each curve $\gamma _{ij}\times \left\{ 1\right\} $ with
framing coefficient -1.  The boundary of $W$ is $\partial W=X_{0}\sqcup -X$
.  We note that
\begin{equation*}
\pi _{1}\left( W\right) =\left\langle x_{1},\ldots ,x_{m}|\left[ x_{i},x_{j}
\right] =1\text{ for all }1\leq i<j\leq m\right\rangle
\cong \mathbb{Z}^{m}\text{.}
\end{equation*}

\begin{figure}[ht!]
\begin{center}
\begin{picture}(0,0)%
\includegraphics{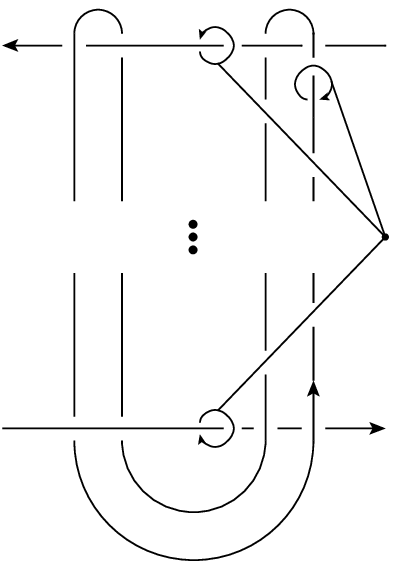}%
\end{picture}%
\setlength{\unitlength}{3947sp}%
\begin{picture}(1871,2661)(-11,-1841)
\small
\put( 30,695){$L_i$}%

\put(1650,-1350){$L_j$}%

\put(800,800){$x_{in_{ij}}$}%

\put(800,-1401){$x_{jn_{ji}}$}%

\put(1666,408){$\overline{\mu}_{ij}$}%

\end{picture}
\nocolon\caption{}
\label{figure3}
\end{center}
\end{figure}

Let $\left\{ x_{ik},\mu _{ijl}\right\} $ be the generators of $\pi
_{1}\left( S^{3}-\left( L\sqcup \gamma \right) \right) $ that are
obtained from a Wirtinger presentation where $x_{ik}$ are
meridians of the $i^{th}$ component of $L$ and $\mu _{ijl}$ are
meridians of the $\left( i,j\right) ^{th}$ component of $\gamma $.
Note that $\left\{ x_{ik},\mu _{ijl}\right\} $ generate $G\equiv
\pi _{1}\left( X\right) $.  For each $ 1\leq i\leq m$ let
$\overline{x}_{i}=x_{i1}$ and $\overline{\mu }_{ij}$ be the
specific $\mu _{ijl}$ that is denoted in Figure \ref{figure3}.
We
will use the
convention that
$$\left[ a,b\right] =aba^{-1}b^{-1}
$$and
$$a^{b}=bab^{-1}\text{.}$$
We can choose a projection of the trivial link so that the arcs $\alpha
_{ij} $ do not pass under a component of $L$.  Since $\overline{\mu }_{ij}$
is equal to a longitude of the curve $\gamma _{ij}$ in $X$, we have
$\overline{\mu }_{ij}=\left[ x_{in_{ij}},\lambda x_{jn_{ji}}\lambda ^{-1}
\right] $ for some $n_{ij}$ and $n_{ji}$ and $\lambda $ where $\lambda $ is
a product of conjugates of meridian curves $\overline{\mu }_{lk}$ and $
\overline{\mu }_{lk}^{-1}$.  Moreover, we can find a projection of $L\sqcup
\gamma $ so that the individual components of $L$ do not pass under or over
one another.  Hence $x_{ij}=\omega \overline{x}_{i}\omega ^{-1}$ where $
\omega $ is a product of conjugates of the meridian curves $\overline{\mu }
_{lk}$ and $\overline{\mu }_{lk}^{-1}$.  As a result, we have
\begin{eqnarray}
\overline{\mu }_{ij} &=&\left[ x_{in_{ij}},\lambda x_{jn_{ji}}\lambda ^{-1}
\right]  \label{muij} \\
&=&\left[ \omega _{1}\overline{x}_{i}\omega _{1}^{-1},\lambda \omega _{2}
\overline{x}_{j}\omega _{2}^{-1}\lambda ^{-1}\right]  \notag \\
&=&\left[ \overline{x}_{i},\omega _{1}^{-1}\lambda \omega _{2}\overline{x}
_{j}\omega _{2}^{-1}\lambda ^{-1}\omega _{1}^{-1}\right] ^{\omega _{1}}
\notag
\end{eqnarray}
for some $\lambda $, $\omega _{1}$, and $\omega _{2}$.

We note that $\overline{\mu }_{ij}=\left[ x_{in_{ij}},\lambda
x_{jn_{ji}}\lambda ^{-1}\right] $ hence $\overline{\mu }_{ij}\in G^{\prime }$
for all $i<j$.  Setting $v=\omega _{1}^{-1}\lambda \omega _{2}$ and using
the equality
\begin{equation}
\left[ a,bc\right] =\left[ a,b\right] \left[ a,c\right] ^b
\end{equation}
we see that
\begin{eqnarray}
\overline{\mu }_{ij} &=&\left[ \overline{x}_{i},v\overline{x}_{j}v^{-1}
\right] ^{\omega _{1}}  \label{muij2} \\
&=&\left[ \overline{x}_{i},v\overline{x}_{j}v^{-1}\right] \text{ mod }
G^{\prime \prime }  \notag \\
&=&\left[ \overline{x}_{i},\left[ v,\overline{x}_{j}\right] \overline{x}_{j}
\right]  \notag \\
&=&\left[ \overline{x}_{i},\left[ v,\overline{x}_{j}\right] \right] \left[
\overline{x}_{i},\overline{x}_{j}\right] ^{\left[ v,\overline{x}_{j}\right]}  \notag \\
&=&\left[ \overline{x}_{i},\left[ v,\overline{x}_{j}\right] \right] \left[
\overline{x}_{i},\overline{x}_{j}\right] \text{ mod }G^{\prime \prime }
\notag
\end{eqnarray}
since $\omega _{1},v\in G^{^{\prime }}$.

Consider the dual relative handlebody decomposition $\left( W,X\right) $.  $
W$ can be obtained from $X$ by adding a 0--framed 2--handle to $X\times I$
along each of the meridian curves $\overline{\mu }_{ij}\times \left\{
1\right\} $.  $\left( \ref{muij}\right) $ implies that $\overline{\mu }
_{ij} $ is trivial in $H_{1}\left( X\right) $ hence the inclusion map $
j\co X\rightarrow W$ induces an isomorphism $j_{\ast }\co H_{1}\left( X\right)
\overset{\cong }{\rightarrow }H_{1}\left( W\right) $.  Therefore if $\phi
\co G\twoheadrightarrow \Lambda $ where $\Lambda $ is abelian then there exists
a $\psi \co \pi _{1}\left( W\right) \twoheadrightarrow \Lambda $ such that $
\psi \circ j_{\ast }=\phi $.

Suppose $\phi \co G\twoheadrightarrow \left\langle t\right\rangle \cong \mathbb{
Z}$ and $\psi \co \pi _{1}\left( W\right) \twoheadrightarrow \left\langle
t\right\rangle $ is an extension of $\phi $ to $\pi _{1}\left( W\right) $. 
Let $X_{\phi }$ and $W_{\psi }$ be the infinite cyclic covers of $W$ and $X$
corresponding to $\psi $ and $\phi $ respectively.  Consider the long exact
sequence of pairs,
\begin{equation}\label{les}
\rightarrow H_{2}\left( W_{\psi },X_{\phi }\right) \overset{\partial _{\ast }
}{\rightarrow }H_{1}\left( X_{\phi }\right) \rightarrow H_{1}\left( W_{\psi
}\right) \rightarrow
\end{equation}
Since $\pi _{1}\left( W\right) \cong \mathbb{Z}^{m}$, $H_{1}\left( W_{\psi
}\right) \cong \mathbb{Z}^{m-1}$ where t acts trivially so that $H_{1}\left( W_{\psi }\right) $
has rank $0$ as
a $\mathbb{Z}\left[ t^{\pm 1}\right] $--module.  $H_{2}\left(
W_{\psi },X_{\phi }\right) \cong \left( \mathbb{Z}\left[ t^{\pm 1}\right]
\right) ^{\binom{m}{2}}$ generated by the core of each 2--handle (extended by
$\overline{\mu }_{ij}\times I$) attached to $X$.  Therefore, $\text{Im}
\partial _{\ast }$ is generated by a lift of $\overline{\mu }_{ij}$ in $
H_{1}\left( X_{\phi }\right) $ for all $1\leq i<j\leq m$.  To show that
$H_{1}\left( X_{\phi }\right)$ has rank $0$
it suffices to show that each
of the $\overline{\mu }_{ij}$ are $\mathbb{Z}\left[ t^{\pm 1}\right] 
$--torsion in $H_{1}\left( X_{\phi }\right) $.

Let $F=\left\langle \overline{x}_{1},\ldots
,\overline{x}_{m}\right\rangle $ be the free group of rank $m$ and
$f\co F\rightarrow G$ be defined by $f\left( \overline{x}_{i}\right)
=\overline{x}_{i}$.  We have the following $\binom{m }{3}$ Jacobi
relations in $F/F^{^{\prime \prime }}$ \cite[Proposition
7.3.6]{kawauchi:knot_theory}.  For all $1\leq i<j<k\leq m$,
\begin{equation*}
\left[ \overline{x}_{i},\left[ \overline{x}_{j},\overline{x}_{k}\right]
\right] \left[ \overline{x}_{j},\left[ \overline{x}_{k},\overline{x}_{i}
\right] \right] \left[ \overline{x}_{k},\left[ \overline{x}_{i},\overline{x}
_{j}\right] \right] =1\text{ mod }F^{\prime \prime }\text{.}
\end{equation*}
Using $f$, we see that these relations hold in $G/G^{\prime \prime }$ as
well.  From $\left( \ref{muij2}\right) $, we can write
\begin{equation*}
\left[ \overline{x}_{i},\overline{x}_{j}\right] =\left[ \left[ v_{ij},
\overline{x}_{j}\right] ,\overline{x}_{i}\right] \overline{\mu }_{ij}\text{
mod }G^{\prime \prime }\text{.}
\end{equation*}
Hence for each $1\leq i<j<k\leq m$ we have the Jacobi relation $J\left(
i,j,k\right) $ in $G/G^{\prime \prime }$,
\begin{eqnarray}
1 &=&\left[ \overline{x}_{i},\left[ \overline{x}_{j},\overline{x}_{k}\right]
\right] \left[ \overline{x}_{j},\left[ \overline{x}_{i},\overline{x}_{k}
\right] ^{-1}\right] \left[ \overline{x}_{k},\left[ \overline{x}_{i},
\overline{x}_{j}\right] \right] \text{ mod }G^{\prime \prime } \notag \\
&=&\left[ \overline{x}_{i},\left[ \left[ v_{jk},\overline{x}_{k}\right] ,
\overline{x}_{j}\right] \overline{\mu }_{jk}\right]
\left[ \overline{x}_{j},
\overline{\mu }_{ik}^{-1}\left[ \overline{x}_{i},\left[ v_{ik},\overline{x}
_{k}\right] \right] \right]
\notag \\ & &
\left[ \overline{x}_{k},\left[ \left[ v_{ij},
\overline{x}_{j}\right] ,\overline{x}_{i}\right] \overline{\mu }_{ij}\right]
\text{ mod }G^{\prime \prime }  \notag \\
&=&\left[ \overline{x}_{i},\left[ \left[ v_{jk},\overline{x}_{k}\right] ,
\overline{x}_{j}\right] \right] \left[ \overline{x}_{i},\overline{\mu }_{jk}
\right] \left[ \overline{x}_{j},\overline{\mu }_{ik}^{-1}\right] \left[
\overline{x}_{j},\left[ \overline{x}_{i},\left[ v_{ik},\overline{x}_{k}
\right] \right] \right]
\notag \\ & &
\left[ \overline{x}_{k},\left[ \left[ v_{ij},
\overline{x}_{j}\right] ,\overline{x}_{i}\right] \right] \left[ \overline{x}
_{k},\overline{\mu }_{ij}\right] \text{ mod }G^{\prime \prime }  \notag \\
&=&\left[ \overline{x}_{i},\overline{\mu }_{jk}\right] \left[ \overline{x}
_{j},\overline{\mu }_{ik}^{-1}\right] \left[ \overline{x}_{k},\overline{\mu }
_{ij}\right] \left[ \overline{x}_{i},\left[ \left[ v_{jk},\overline{x}_{k}
\right] ,\overline{x}_{j}\right] \right] \left[ \overline{x}_{j},\left[
\overline{x}_{i},\left[ v_{ik},\overline{x}_{k}\right] \right] \right]
\notag \\ & &
\left[\overline{x}_{k},\left[ \left[ v_{ij},\overline{x}_{j}\right] ,\overline{x}
_{i}\right] \right] \text{ mod }G^{\prime \prime }\text{.}  \label{freerel}
\end{eqnarray}

\begin{figure}[ht!]
\begin{center}
\begin{picture}(0,0)%
\includegraphics{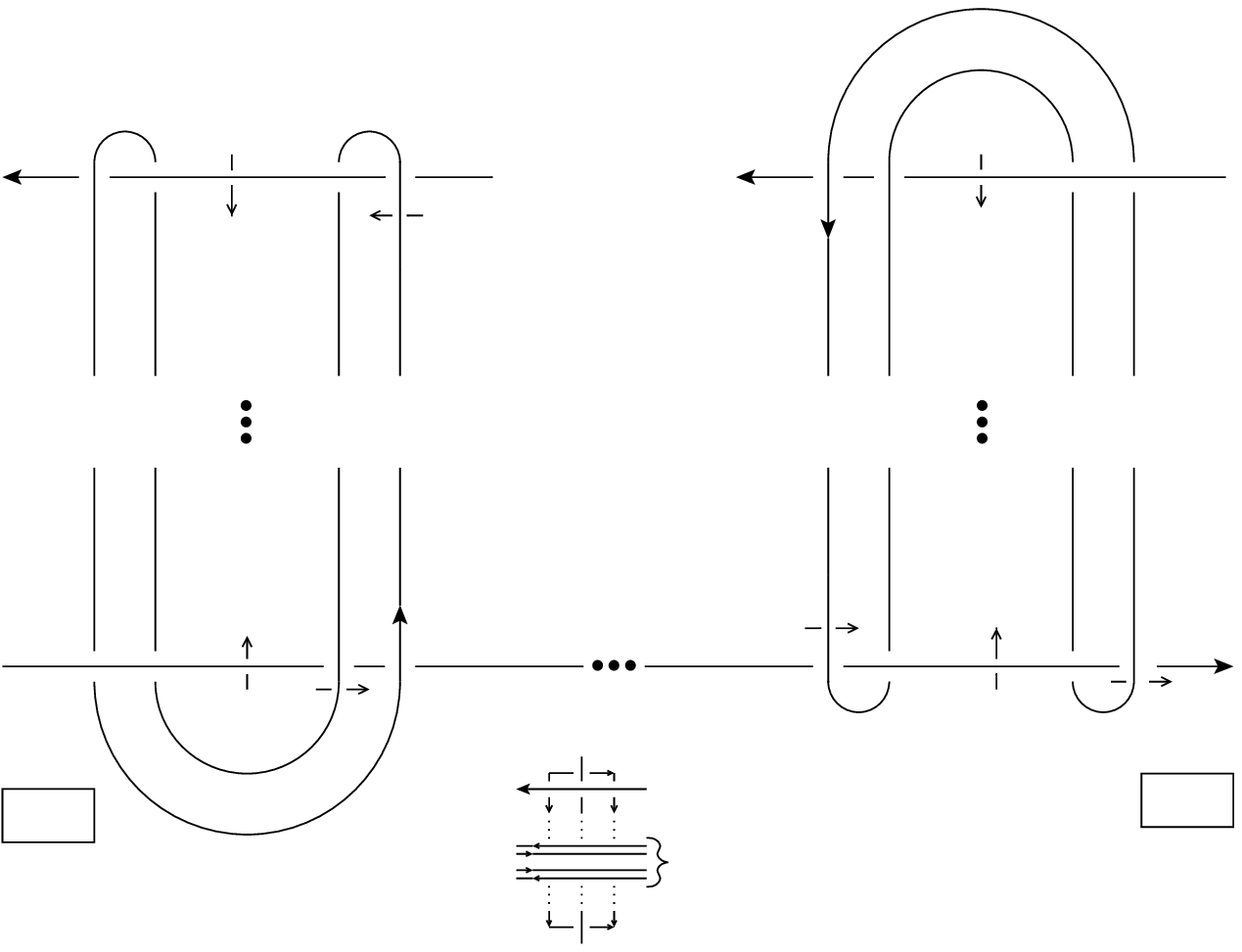}%
\end{picture}%
\setlength{\unitlength}{3947sp}%
\begin{picture}(6061,4670)(-49,-3858)

\put(2850,-2880){\makebox(0,0)[lb]{\smash{\SetFigFont{6}{6}{\familydefault}{\mddefault}{\updefault}{$\overline{\mu}_{ji}$}%
}}}
\put(3250,-3400){\makebox(0,0)[lb]{\smash{\SetFigFont{6}{6}{\familydefault}{\mddefault}{\updefault}{$\gamma_j$}%
}}}
\put(3003,-3110){\makebox(0,0)[lb]{\smash{\SetFigFont{6}{6}{\familydefault}{\mddefault}{\updefault}{$x_{jn_{ji}}$}%
}}}
\put(2604,-3775){\makebox(0,0)[lb]{\smash{\SetFigFont{6}{6}{\familydefault}{\mddefault}{\updefault}{$\alpha_j$}%
}}}
\put(2326,-2986){\makebox(0,0)[lb]{\smash{\SetFigFont{6}{6}{\familydefault}{\mddefault}{\updefault}{$L_j$}%
}}}

\small
\put(2513,-2350){$L_i$}%

\put( 76, 60){$L_j$}%

\put(910, 150){$x_{jn_{ji}}$}%

\put(2063,-286){$\overline{\mu}_{ji}$}%

\put(938,-2236){$x_{in_{ij}}$}%

\put(1314,-2625){$\alpha_j$}%


\put(-20,-3211){$j<i$}%

\put(3620,-2236){$\overline{\mu}_{ik}$}%

\put(4687,-2660){$x_{in_{ik}}$}%

\put(5700,-2611){$\beta$}%

\put(5550,-3136){$k>i$}%

\put(4575,150){$x_{kn_{ki}}$}%

\put(3637, 50){$L_k$}%

\end{picture}
\nocolon\caption{}
\label{figure4}
\end{center}
\end{figure}

Moreover, for each component of the trivial link $L_{i}$ the
longitude, $ l_{i}$, of $L_{i}$ is trivial in $G$ and is a product
of commutators of $ \overline{\mu }_{ij}$ with a conjugate of
$\overline{x}_{j}$. We can write each of the longitudes (see
Figure \ref{figure4})
as
\begin{eqnarray}
l_{i} &=&\prod_{j<i}\alpha _{j}\lambda _{j}^{-1}\overline{\mu }
_{ji}^{-1}\lambda _{j}\cdot \prod_{k>i}\overline{\mu }_{ik}\beta _{k}\text{
mod }G^{\prime \prime } \notag \\
&=&\prod_{j<i}\left( \lambda _{j}^{-1}x_{jn_{ji}}^{-1}\overline{\mu }
_{ji}x_{jn_{ji}}\lambda _{j}\right) \lambda _{j}^{-1}\overline{\mu }
_{ji}^{-1}\lambda _{j}\cdot \notag \\
& & \prod_{k>i}\overline{\mu }_{ik}\left( \lambda
_{k}x_{kn_{ki}}^{-1}\lambda _{k}^{-1}\overline{\mu }_{ik}^{-1}\lambda
_{k}x_{kn_{ki}}\lambda _{k}^{-1}\right) \notag \\
&=&\prod_{j<i}\left[ x_{jn_{ji}}^{-1},\overline{\mu }_{ji}\right] ^{\lambda
_{j}^{-1}}\cdot \prod_{k>i}\left[ \overline{\mu }_{ik},\lambda
_{k}x_{kn_{ki}}^{-1}\lambda _{k}^{-1}\right] \notag \\
&=&\prod_{j<i}\left[ \overline{x}_{j}^{-1},\overline{\mu }_{ji}\right] \cdot
\prod_{k>i}\left[ \overline{\mu }_{ik},\overline{x}_{k}^{-1}\right] \text{
mod }G^{\prime \prime }\text{.} \label{longrel}
\end{eqnarray}

It follows that
\begin{equation*}
\prod_{j<i}\left[ \overline{x}_{j}^{-1},\overline{\mu }_{ji}\right] \cdot
\prod_{k>i}\left[ \overline{\mu }_{ik},\overline{x}_{k}^{-1}\right] =1\text{
mod }G^{\prime \prime }\text{.}
\end{equation*}

Since $G^{\prime \prime }\subset \left[ \ker \phi ,\ker \phi \right] $, the
relations in (\ref{freerel}) and (\ref{longrel}) hold in $H_{1}\left( X_{\phi
}\right)$\nl$\left( = \ker \phi /\left[ \ker \phi ,\ker \phi \right] \right)$ as well. 
Suppose $\phi \co G\twoheadrightarrow \mathbb{Z}$ is defined by sending $\overline{x}
_{i}\longmapsto t^{n_{i}}$.  Since $\phi $ is surjective, $n_{N}\neq 0$
for
some $N$.
We consider a subset of $\binom{m}{2}$ relations in $H_{1}\left( X_{\phi
}\right) $ that we index by $\left( i,j\right) $ for $1\leq i<j\leq m$. 
When $i=N$ or $j=N$ we consider the $m-1$ relations
$$\text{(i)}\qquad R_{iN}=l_{i}
\qquad
\text{and}
\qquad
\text{(ii)}\qquad R_{Nj}=l_{j}^{-1}\text{.}$$
Rewriting $l_{i}$ as an element of the $\mathbb{Z}\left[ t^{\pm 1}\right] $
-module $H_{1}\left( X_{\phi }\right) $ generated by \nl$\left\{ \overline{\mu
}
_{ij}|1\leq i<j\leq m\right\}$ from (\ref{longrel}) we have
\begin{eqnarray}
R_{iN} &=&\sum_{j<i}\left( t^{-n_{j}}-1\right) \overline{\mu }
_{ji}+\sum_{k>i}\left( 1-t^{-n_{k}}\right) \overline{\mu }_{ik} \notag
\\
&=&\sum_{j<i}t^{-n_{j}}\left( 1-t^{n_{j}}\right) \overline{\mu }
_{ji}+\sum_{k>i}t^{-n_{k}}\left( t^{n_{k}}-1\right) \overline{\mu }_{ik}
\notag \\
&=&\sum_{j<i}\left[ \left( 1-t^{n_{j}}\right) +\left( t^{-n_{j}}-1\right)
\left(1- t^{n_{j}}\right) \right] \overline{\mu }_{ji}+ \notag \\
& &
\sum_{k>i}\left[
\left( t^{n_{k}}-1\right) +\left( t^{-n_{k}}-1\right) \left(
t^{n_{k}}-1\right) \right] \overline{\mu }_{ik}\text{.}  \label{r1}
\end{eqnarray}
Similarily, we have
\begin{eqnarray}
R_{Nj}&=&\sum_{i<j}\left[ \left( t^{n_{i}}-1\right) +\left(
t^{-n_{i}}-1\right) \left( t^{n_{i}}-1\right) \right] \overline{\mu }
_{ij}+ \notag \\
& &
\sum_{k>j}\left[ \left( 1-t^{n_{k}}\right) +\left( t^{-n_{k}}-1\right)
\left( 1-t^{n_{k}}\right) \right] \overline{\mu }_{jk}\text{.}  \label{r2}
\end{eqnarray}
For the other $\binom{m-1}{3}$ relations, we use the Jacobi relations from $
\left( \ref{freerel}\right) $.  Define $R_{ij}$ to be
\begin{equation*}
R_{ij}=\left\{
\begin{array}{c}
J\left( N,i,j\right) \text{ for }N<i<j \\
J\left( i,N,j\right) ^{-1}\text{ for }i<N<j \\
J\left( i,j,N\right) \text{ for }i<j<N
\end{array}
\text{.}\right.
\end{equation*}
We can write these relations as
\begin{equation}
R_{ij}=\left\{
\begin{array}{l}
\left( t^{n_{j}}-1\right) \overline{\mu }_{Ni}+\left( 1-t^{n_{i}}\right)
\overline{\mu }_{Nj}+\left( t^{n_{N}}-1\right) \overline{\mu }_{ij} + \\
\left( t^{n_{N}}-1\right) \left( t^{n_{i}}-1\right) \left(
t^{n_{j}}-1\right) \left(
\widetilde{v}_{ij}+\widetilde{v}_{Nj}-\widetilde{v}
_{Nj}\right) \text{ for }N<i<j \\
\left( 1-t^{n_{j}}\right) \overline{\mu }_{iN}+\left( t^{n_{N}}-1\right)
\overline{\mu }_{ij}+\left( 1-t^{n_{i}}\right) \overline{\mu }_{Nj} + \\
\left( t^{n_{N}}-1\right) \left( t^{n_{i}}-1\right) \left(
t^{n_{j}}-1\right) \left(-
\widetilde{v}_{iN}-\widetilde{v}_{Nj}+\widetilde{v}
_{ij}\right) \text{ for }i<N<j \\
\left( t^{n_{N}}-1\right) \overline{\mu }_{ij}+\left( 1-t^{n_{j}}\right)
\overline{\mu }_{iN}+\left( t^{n_{i}}-1\right) \overline{\mu }_{jN} + \\
\left( t^{n_{N}}-1\right) \left( t^{n_{i}}-1\right) \left(
t^{n_{j}}-1\right) \left(
\widetilde{v}_{ij}+\widetilde{v}_{jN}-\widetilde{v}
_{iN}\right) \text{ for }i<j<N
\end{array}
\right.  \label{r3}
\end{equation}
\newline
where $\widetilde{v}_{ij}$ is a lift of $v_{ij}$.

For $1\leq i<j\leq m$ order the pairs $ij $ by the dictionary
ordering. That is, $ij<lk$ provided either $i<l$ or $j<k$ when $i=l$.
The relations above give us an $\binom{m}{2}\times \binom{m}{2}$
matrix $M$ with coefficients in $\mathbb{Z}\left[ t^{\pm 1}\right] $.  The $
\left( ij,kl\right) ^{th}$ component of $M$ is the coefficient of $\overline{
\mu }_{kl}$ in $R_{ij}$.  We claim for now that
\begin{equation}
M=\left( t^{n_{N}}-1\right) I+\left( t-1\right) S+\left( t-1\right) ^{2}E
\label{matrixform}
\end{equation}
for some ``error'' matrix $E$ where $I$ is the identity matrix and $S$ is a
skew-symmetric matrix.  For an example, when $m=4$ and $N=1$, $M$ is the
matrix
\begin{equation*}
\left[
\begin{array}{cccccc}
t^{n_{1}}-1 & 0 & 0 & 1-t^{n_{3}} & 1-t^{n_{4}} & 0 \\
0 & t^{n_{1}}-1 & 0 & t^{n_{2}}-1 & 0 & 1-t^{n_{4}} \\
0 & 0 & t^{n_{1}}-1 & 0 & t^{n_{2}}-1 & t^{n_{3}}-1 \\
t^{n_{3}}-1 & 1-t^{n_{2}} & 0 & t^{n_{1}}-1 & 0 & 0 \\
t^{n_{4}}-1 & 0 & 1-t^{n_{2}} & 0 & t^{n_{1}}-1 & 0 \\
0 & t^{n_{4}}-1 & 1-t^{n_{3}} & 0 & 0 & t^{n_{1}}-1
\end{array}
\right]
+\left( t-1\right) ^{2}E\text{.}
\end{equation*}
The proof of $\left( \ref{matrixform}\right) $ is left until the end.

We will show that $M$ is non-singular as a matrix over the quotient field $
\mathbb{Q}\left( t\right) $.  Consider the matrix $A=\frac{1}{t-1}M$.  We
note that $A$ is a matrix with entries in $\mathbb{Z}\left[ t^{\pm
1}\right]
$ and $A\left( 1\right) $ evaluated at $t=1$ is
\begin{equation*}
A\left( 1\right) =NI+S\left( 1\right) \text{.}
\end{equation*}
To show that $M$ is non-singular, it suffices to show that $A\left( 1\right)
$ is non-singular.

Consider the quadratic form $q\co \mathbb{Q}^{\binom{m}{2}}\rightarrow \mathbb{Q
}^{\binom{m}{2}}$ defined by $q\left( z\right) \equiv z^{T}A\left( 1\right) z
$ where $z^{T}$ is the transpose of $z$.  Since $A\left( 1\right) =N
I+S\left( 1\right) $ where $S\left( 1\right) $ is skew-symmetric we have,
\begin{equation*}
q\left( z\right) =N\sum z_{i}^{2}\text{.}
\end{equation*}
Moreover, $N \neq 0$ so $q\left( z\right) =0$ if and only if
$z=0$.  Let $z$ be a vector
satisfying $A\left( 1\right) z=0$.  We have $q\left( z\right) =z^{T}A\left(
1\right) z=z^{T}0=0$ which implies that $z=0$.  Therefore $M$ is a
non-singular matrix.  This implies that each element $\overline{\mu }_{ij}$
is $\mathbb{Z}\left[ t^{\pm 1}\right] $--torsion which will complete the the
proof once we have established the above claim.

We ignore entries in $M$ that lie in $J^{2}$ where $J$ is the augmentation
ideal of $\mathbb{Z}\left[ t^{\pm 1}\right] $ since they only contribute
to
the error matrix $E$.  Using $\left( \ref{r1}\right) $, $\left( \ref{r2}
\right) $, and $\left( \ref{r3}\right) $ above we can explicitely write
the
entries in $M$ $\left( \text{mod }J^{2}\right) $.  Let $m_{ij,lk}$
denote the $\left( ij,lk\right) $ entry of $M$ $\left( \text{mod }
J^{2}\right) $. 

Case 1 $\left( j=N\right) $:\qua From $\left( \ref{r1}\right) $ we have
\begin{equation*}
m_{iN,li}=1-t^{n_{l}}\text{, }m_{iN,ik}=t^{n_{k}}-1\text{,}
\end{equation*}
and $m_{iN,lk}=0$ when neither $l$ nor $k$ is equal to $N$.

Case 2 $\left( i=N\right) $:\qua  From $\left( \ref{r2}\right) $ we have
\begin{equation*}
m_{Nj,lj}=t^{n_{l}}-1\text{, }m_{Nj,jk}=1-t^{n_{k}}\text{,}
\end{equation*}
and $m_{Nj,lk}=0$ when neither $l$ nor $k$ is equal to $N$.

Case 3 $\left( N<i<j\right) $:\qua  From $\left( \ref{r3}\right) $ we have
\begin{equation*}
m_{ij,Ni}=t^{n_{j}}-1\text{, }m_{ij,Nj}=1-t^{n_{i}}\text{, }
m_{ij,ij}=t^{n_{N}}-1\text{,}
\end{equation*}
and $m_{ij,lk}=0$ otherwise.

Case 4 $\left( i<N<j\right) $:\qua  From $\left( \ref{r3}\right) $ we have
\begin{equation*}
m_{ij,iN}=1-t^{n_{j}}\text{, }m_{ij,ij}=t^{n_{N}}-1\text{, }
m_{ij,Nj}=1-t^{n_{i}}\text{,}
\end{equation*}
and $m_{ij,lk}=0$ otherwise.

Case 5 $\left( i<j<N\right) $:\qua  From $\left( \ref{r3}\right) $ we have
\begin{equation*}
m_{ij,ij}=t^{n_{N}}-1\text{, }m_{ij,iN}=1-t^{n_{j}}\text{, }
m_{ij,jN}=t^{n_{i}}-1\text{,}
\end{equation*}
and $m_{ij,lk}=0$ otherwise. 

We first note that in each of the cases, the diagonal entries $m_{ij,ij}$
are all $t^{n_{N}}-1$.  Next, we will show that the off diagonal entries
have the property that $m_{ij,lk}=-m_{lk,ij}$ for $ij<lk$.  This will
complete the proof of the claim since we see that each entry is divisible
by $t-1$.

We verify the skew symmetry in Cases 1 and 3.  The other cases are
similar
and we leave the verifications to the reader.

Case 1 $\left( j=N\right) $:
\begin{equation*}
m_{iN,li}=1-t^{n_{l}}=-m_{li,iN}\text{ (case 5)}
\end{equation*}
and
\begin{equation*}
m_{iN,ik}=t^{n_{k}}-1=-m_{ik,iN}\text{ (case 4).}
\end{equation*}

Case 3 $\left( N<i<j\right) $:
\begin{equation*}
m_{ij,Ni}=t^{n_{j}}-1=-m_{Ni,ij}\text{ (case 2)}
\end{equation*}
and
$$m_{ij,Nj}=1-t^{n_{i}}=-m_{Nj,ij}\text{ (case 2).}\eqno{\qed}$$

\begin{proposition}\label{lcs}Let $X$ be as in Theorem \ref{mainthm},
$G=\pi_1\left(X\right)$ and $F$ be the
free group on $2$ generators.  There is no epimorphism from $G$ onto $F/F_4$. \end{proposition}

\begin{proof}
Let $F=\left\langle x,y\right\rangle $ be the free group and $\phi
\co F\left/
F_{4}\right. \twoheadrightarrow \left\langle t\right\rangle $ be defined
by $
x\longmapsto t$ and $y\longmapsto 1$. Suppose that there exists a
surjective map $\eta \co G\twoheadrightarrow F\left/ F_{4}\right. $.  Let
$N=\ker \phi $ and $H=\ker \left(\eta \circ \phi \right)$.  Since $\eta $
is
surjective we get an epimorphism of $\mathbb{Z}\left[ t^{\pm 1}
\right] $--modules $\widetilde{\eta }\co H\left/ H^{\prime }\right.
\twoheadrightarrow N\left/ N^{\prime }\right. $. From $\left( \ref{les}
\right) $ we get the short exact sequence
\begin{equation*}
0\rightarrow \text{Im}\partial _{\ast }\overset{i}{\rightarrow
}H_{1}\left(
X_{\eta \circ \phi }\right) \rightarrow H_{1}\left( W_{\psi }\right)
\rightarrow 0.
\end{equation*}
Let $J$ be the augmentation ideal of $\mathbb{Z}\left[ t^{\pm 1}\right]
$.
We compute $N\left/ N^{\prime }\right. \cong \mathbb{Z}\left[ t^{\pm 1}%
\right] \left/ J^{3}\right. $ so that $\widetilde{\eta }\co H_{1}\left(
X_{\eta
\circ \phi }\right) \twoheadrightarrow \mathbb{Z}\left[ t^{\pm 1}\right]
\left/ J^{3}\right. $.  Let $\sigma \in H_{1}\left( X_{\eta \circ \phi
}\right) $ such that $\widetilde{\eta }\left( \sigma \right) =1$. Since
every element in $H_{1}\left( W_{\psi }\right) \cong
\bigoplus_{i=1}^{m-1}%
\frac{\mathbb{Z}\left[ t^{\pm 1}\right] }{J}$ is $\left( t-1\right) $%
--torsion, $\left( t-1\right) \sigma \in \text{Im}\partial _{\ast }$ hence
$t-1\in \text{Im}\left(\widetilde{\eta }\circ i\right)$. Recall that in
the proof of
the
Theorem \ref{mainthm}, we showed that there exists a surjective
$\mathbb{Z}%
\left[ t^{\pm 1}\right] $--module homomorphism $\rho \co P\twoheadrightarrow
\text{Im}\partial _{\ast }$ $\ $where $P$ is finitely presented as
\begin{equation*}
0\rightarrow \mathbb{Z}\left[ t^{\pm 1}\right] ^{\binom{m}{2}}\overset{%
\left( t-1\right) A}{\rightarrow }\mathbb{Z}\left[ t^{\pm 1}\right]
^{\binom{%
m}{2}}\overset{\pi }{\rightarrow }P\rightarrow 0\text{.}
\end{equation*}
Let $g\co P\rightarrow \mathbb{Z}\left[ t^{\pm 1}\right] \left/ J^{3}\right.
$
defined by $g\equiv \widetilde{\eta }\circ i\circ \rho $. Since $\rho $
is
surjective, $t-1\in \text{Im}g$. After tensoring with $\mathbb{Q}\left[
t^{\pm 1}\right] $, we get a map $g\co P\otimes _{\mathbb{Z}\left[ t^{\pm 1}%
\right] }\mathbb{Q}\left[ t^{\pm 1}\right] \rightarrow \mathbb{Q}\left[
t^{\pm 1}\right] \left/ J^{3}\right. $. It is easy to see that either
$g$
is surjective or the image of $g$ is the submodule generated by $t-1$. 
Note that the submodule generated by $t-1$ is isomorphic $\mathbb{Q}\left[
t^{\pm 1}\right] \left/ J^{2}\right. $. Hence, in either case, we get a
surjective map $h\co P\otimes _{\mathbb{Z}\left[ t^{\pm 1}\right]
}\mathbb{Q}%
\left[ t^{\pm 1}\right] \rightarrow \mathbb{Q}\left[ t^{\pm 1}\right]
\left/
J^{2}\right. $.

Consider the $\mathbb{Q}\left[ t^{\pm 1}\right] $--module $P^{\prime }$
presented by $A$. Let $h^{\prime }\co \mathbb{Q}\left[ t^{\pm 1}\right] ^{%
\binom{m}{2}}\rightarrow$\break $\mathbb{Q}\left[ t^{\pm 1}\right] \left/
J^{2}\right. $ be defined by $h^{\prime }=\left( t-1\right) h\circ \pi $.
Since
\begin{equation*}h^{\prime }\left( A\left( \sigma \right) \right) =\left(
t-1\right)
h\left( \pi \left( A\left( \sigma \right) \right) \right) =h\left( \pi
\left( \left( t-1\right) A\left( \sigma \right) \right) \right) =
h\left(0\right) =0,\end{equation*} this defines a map $h^{\prime
}\co P^{\prime
}\rightarrow \mathbb{Q}\left[ t^{\pm 1}\right] \left/ J^{2}\right. $ whose
image is the submodule generated by $t-1$. It follows that $P^{\prime}$
maps onto $\mathbb{Q}\left[ t^{\pm 1}\right] \left/ J\right. $.
Setting $t=1$, the vector space over $\mathbb{Q}$ presented by
$A\left(1\right)$ maps onto $\mathbb{Q}$.
Therefore $det(A(1))=0$.
However, it was previously shown that $A\left(1\right)$ was non-singular
which is a
contradiction.
\end{proof}

\begin{corollary}For any closed, orientable 3--manifold $Y$ with $P/P_{4}
\cong G/G_{4}$ where $P=\pi _{1}\left( Y\right)$ and $G=\pi _{1}\left(
X\right)$ is the fundamental group of the examples in Theorem
\ref{mainthm}, $c(Y)=1$.
\end{corollary}

Using Proposition \ref{lcs}, it is much easier to show that there exist
\emph{hyperbolic} 3--manifolds with cut number 1.
\begin{corollary}
For each $m\geq 1$ there exist closed, orientable, hyperbolic 3--manifolds
$Y$ with $\beta _{1}\left( Y\right) =m$ such that $\pi_1\left(Y\right)$
cannot map onto $F/F_{4}$ where $F$ is the free group on 2 generators.
\end{corollary}

\begin{proof}
Let $X$ be one of the 3--manifolds in Theorem \ref{mainthm}.  By
\cite[Theorem 2.6]{ruberman:seifert}, there exists a degree one
map $f\co Y\rightarrow X$ where $Y$ is hyperbolic and $f_{\ast }$ is
an isomorphism on $H_{\ast }$. Denote by $G=\pi _{1}\left(
X\right) $ and $P=\pi _{1}\left( Y\right) $. It follows from
Stalling's theorem \cite{stallings:centralseries} that $f$ induces
an isomorphism $f_{*}\co P/P_{n} \rightarrow G/G_{n}$.  In particular
this is true for $n=4$ which completes the proof.
\end{proof}

\end{document}